\documentclass[leqno,12pt]{article}
\usepackage{amsmath, amssymb}
\usepackage{amsthm}

\theoremstyle{plain}
\newtheorem{theorem}{\indent\sc Theorem}[section]
\newtheorem{lemma}[theorem]{\indent\sc Lemma}
\newtheorem{corollary}[theorem]{\indent\sc Corollary}
\newtheorem{proposition}[theorem]{\indent\sc Proposition}

\theoremstyle{definition}
\newtheorem{definition}[theorem]{\indent\sc Definition}
\newtheorem{remark}[theorem]{\indent\sc Remark}
\newtheorem{example}[theorem]{\indent\sc Example}

\makeatletter

\def\address#1#2{\begingroup
\noindent\parbox[t]{7.8cm}{%
\small{\scshape\ignorespaces#1}\par\vskip1ex
\noindent\small{\itshape E-mail address}%
\/: #2\par\vskip4ex}\hfill%
\endgroup}%
\makeatother
\title{\uppercase{The main component of the toric Hilbert scheme}}
\author{ \bigskip \\ \textsc{Olga V. Chuvashova$^{*}$}}

\date{}

\newcommand{\T}{\mathbb{T}}

\newcommand{\X}{\mathbb{X}}
\newcommand{\W}{\mathbb{W}}
\newcommand{\Sp}{{\rm Spec\,}}
\newcommand{\PP}{{\rm Proj\,}}
\newcommand{\R}{\mathbb{R}}
\newcommand{\Z}{\mathbb{Z}}
\newcommand{\N}{\mathbb{N}}
\newcommand{\A}{\mathbb{A}}
\newcommand{\G}{\mathbb{G}_m}
\newcommand{\Chi}{\mathcal{X}}
\renewcommand{\Pr}{\mathbb{P}}
\renewcommand{\l}{\lambda}

\renewcommand{\c}{\cdot}

\newcommand{\Hom}{{\rm Hom}}

\newcommand{\C}{\mathcal{C}}

\renewcommand{\O}{\mathcal{O}}
\newcommand{\cone}{\rm{cone}}

\newcommand{\la}{\langle}
\newcommand{\ra}{\rangle}

\begin{document}

\maketitle

\footnote{$^{*}$ The work was carried out while the author was visiting the Institute Fourier (Grenoble) and supported by grant of the Research Training Network LIEGRITS.}

\begin{abstract}

Let $\X$ be an affine toric variety with  big torus $\T\subset \X$ and let $T\subset\T$ be a subtorus. The general
$T$-orbit closures in $\X$ and their flat limits are parametrized by the main component $H_0$ of the toric
Hilbert scheme. Further, the quotient torus
$\T/T$ acts on $H_0$ with a dense orbit. We describe the fan of this toric variety; this leads us to an integral
analogue of the fiber polytope of Billera and Sturmfels. We also describe the relation of $H_0$ to the main component of the inverse limit of GIT quotients of $\X$ by $T$.

\end{abstract}

\section{Introduction}

The multigraded Hilbert scheme parametrizes, in a technical sense specified below, all homogeneous  ideals in a
polynomial algebra (or, more generally, in an arbitrary finitely generated algebra)  having a fixed Hilbert function
with respect to a grading by an abelian group. In \cite{HS} it was shown that the multigraded Hilbert scheme
always exists as a quasiprojective scheme.

We consider the following case. Let $\X$
be an affine toric (not necessarily normal) variety with big torus $\T\subset \X$ and let
$T\subset \T$ be a subtorus acting on $\X$ by the restriction of
the action of $\T$. This  defines a grading of the algebra of
regular functions $k[\X]$ by the group  of
characters of $T$. Denote by $H_{\X, T}$ the toric Hilbert scheme, that
is, the multigraded Hilbert scheme parametrizing those $T$-invariant
ideals in $k[\X]$ having the same Hilbert function as the toric
$T$-variety $X=\overline{Tx}$, where $x\in \X$ lies in the open
$\T$-orbit \cite{PS}. There is a canonical irreducible component
$H_0$ of $H_{\X, T}$ parametrizing general $T$-orbit closures in
$\X$ and their flat limits (Proposition \ref{prop1}(2)). This component contains
an open orbit for a natural action of $\T/T$ on $H_{\X, T}$. The
main result of this work is a description of the fan of this toric variety (Theorem \ref{th1}). Also, we compare the fan of $H_0$ with
the fan of the toric Chow quotient.

The Chow quotient of a
 {\it projective} toric variety was considered in \cite{KSZ}. In
 particular, in this paper there is a description of its fan. Namely, recall
 that the fan of a projective toric  variety is the normal fan of a
 convex polytope $P$ in the space  generated
 by the lattice of characters of $\T$. Let $Q$ be the projection of
 this polytope on the space $\Chi(T)_{\mathbb{R}}$ generated by
 the lattice of characters of the subtorus $T$. Then the fan of
 the Chow quotient is the normal fan to the fiber polytope $F(P, Q)$ \cite{BS}, which , in a well-defined sense,
 is the average over all fibers of the projection of $P$ on $Q$.
More generally, the fiber fan for a projection of an arbitrary
polyhedron was defined in  \cite{CM}. In this paper some results
of \cite{KSZ} were generalized on the case of a variety that is
projective over some affine variety.

In our affine setting, we show that the fan corresponding to the toric variety $H_0$ is
the normal fan to the average over all "integral" fibers of the
corresponding cone projection. Here by an integral fiber  we mean
the polyhedron generated by all integral points of a fiber of the projection.
Thus this object can be regarded as an integral analogue of the
fiber fan. If $\X$ is a finite-dimensional $T$-module and the grading of $k[\X]$ by the weights of $T$ is positive, then  the fan of $H_0$ coincides with the
normal fan to the state polytope of Sturmfels (see \cite[Theorem 2.5]{St}).

In the last section we consider the toric Chow morphism from the Hilbert scheme to the inverse limit of GIT
quotients $\X/_{\!\!\chi} T$. This morphism was constructed in \cite[Section 5]{HS} in the case when $\X$ is a finite-dimensional $T$-module. We generalize this to the case of a normal affine toric $\T$-variety $\X$ (Theorem \ref{chow}).

The author is indebted to Michel Brion for posing the problem,
numerous ideas and helpful remarks. Thanks are also due to Ivan V.
Arzhantsev for permanent support and useful discussions. The
author is grateful to the referee for important comments which
leaded to improvement of the article.

\section{Terminology and notation}

We consider the category of schemes over an algebraically closed field $k$. A {\it variety} is a separated
integral scheme of finite type. Recall that any scheme $Z$ is characterized by its {\it functor of points} from the category of
$k$-algebras to the category of sets:
$$\underline{Z}: \underline{k-Alg}\to \underline{Set}, \ \ \
\underline{Z}(R):={\rm Mor}(\Sp R, Z),$$ where ${\rm Mor}(\Sp R, Z)$ is the set of morphisms of schemes over $k$
from $\Sp R$ to $Z$ (we denote the functor
 of points of a scheme by the corresponding underlined letter).
Our main reference on schemes is \cite{EH}.
  We denote by $\O_Z$ the structure sheaf of $Z$, and if $Z$ is affine,
then $k[Z]$ denotes the algebra of sections of $\O_Z$ over $Z$.
 We denote by $\A^n$ the affine space $\Sp k[x_1,\ldots, x_n]$.

An $n$-dimensional {\it torus} $T$ is an algebraic group isomorphic to the direct product of $n$ copies of
the multiplicative group $\G$ of the field $k$.
 For the lattices of characters and one-parameter subgroups of $T$, we use the
 notations $\Chi(T) = \Hom(T, \G)$ and $\Lambda(T)= \Hom(\G, T)$.
 We denote by $\la \cdot, \cdot \ra$  the natural pairing between $\Chi(T)$ and $\Lambda(T)$.
For a lattice $\Chi$, let $\Chi_{\R}=\Chi\otimes_{\Z}\R$. If $\Sigma \subset \Chi$ is a monoid, then
 $\cone(\Sigma)$ denotes the cone in
$\Chi_{\R}$ generated by $\Sigma$. For subsets $D_1, D_2$ of a vector space,  we denote by $D_1+D_2$ the
Minkowski sum.

By a {\it toric variety} under a torus $T$  we mean a variety $X$ such that  $T$ is
 embedded as an open subset into $X$,
 the action of $T$  on itself by multiplication extends to an action on $X$, and $X$ admits an open covering by affine $T$-invariant charts.
  We do not require $X$ to be normal.

We denote by $\C_X$ the associated fan of a toric variety $X$, so the cones of $\C_X$ lie in $\Lambda(T)_{\R}$ (see \cite[Sec.
1.4]{Ful}). The $T$-orbits on $X$ are in order-reversing one-to-one correspondence with the cones of $\C_X$. If
$\sigma(Y)$ is the cone in $\C_X$ corresponding to a $T$-orbit  $Y$, then a one-parameter subgroup $\lambda
\in \Lambda(T)$ lies in the interior of $\sigma(Y)$ if and only if $\lim_{s\to 0}\lambda(s)$ exists and
 lies in $Y$. A toric variety is
 determined by its fan up to normalization.

\section{Definitions
and background on multigraded Hilbert schemes}

Let $\X$ be an affine variety over $k$ with an action of a
torus $T$, so its algebra of regular functions $S:=k[\X]$ is graded by the group $\Chi(T)$ of characters of
$T$ :
$$S=\bigoplus_{\chi\in\, \Chi(T)} S_{\chi},$$ where  $S_{\chi}$ is the
subspace of $T$-semiinvariant functions of weight $\chi$. Let
$$\Sigma:=\{\chi\in \Chi(T)\ ;\ S_\chi\ne 0\}.$$
This is a finitely generated monoid. Conversely, if $S$ is a finitely generated commutative $k$-algebra without zero divizors
graded by $\Chi(T)$, then we have a $T$-action on the affine  variety $\X=\Sp S$.

\begin{definition}
The grading of $S$ by $\Chi(T)$ is {\it positive} if $k[\X]_0=k$ and  $\cone(\Sigma)$ is strictly convex.
\end{definition}

Notice that in
the case of a positive grading there exists a unique minimal system of generators of $\Sigma$.
The following definition was introduced in \cite{HS}.

\begin{definition}\label{Def}  Given a function $h:\Chi(T)\to \N$, the {\it Hilbert functor}
is the covariant functor $\underline{H^h_{\X, T}}$ from the category of $k$-algebras to the category of sets
assigning to any $k$-algebra $R$ the set of all $T$-invariant ideals $I\subseteq R\otimes_k S$ such that
$(R\otimes_k S_\chi)/I_\chi$ is a locally free $R$-module of rank $h(\chi)$ for any $\chi\in \Chi(T)$.
\end{definition}

Remark that we can also view $\underline{H^h_{\X,
T}}(R)$ as a set of closed $T$-invariant subschemes $Y\subset \Sp R \times \X$ such that the projection $Y\to \Sp R$ is flat.

In \cite[Theorem 1.1]{HS} it was proved that there exists a quasiprojective scheme $H^h_{\X, T}$ which represents
this functor in the case when $\X$ is a finite-dimensional $T$-module $V$.
In the case of an arbitrary $\X$ there
exists a $T$-equivariant closed immersion $\X\hookrightarrow V$, where $V$ is a finite-dimensional $T$-module.
Then the Hilbert functor $\underline{H^h_{\X, T}}$ is represented by a closed subscheme of $H^h_{V, T}$ (see
\cite[Lemma 1.6]{AB}). Namely, for an algebra $R$ the subset
$\underline{H^h_{\X, T}}(R)\subset \underline{H^h_{V, T}}(R)$ consists of those ideals $I\subset R\otimes_k k[V]$
that $I\in \underline{H^h_{V, T}}(R)$ and $R\otimes_k I_{\X}\subset I$.

Recall that the {\it universal family} is the closed subscheme $\W_{\X, T}$ of $H_{\X, T}\times \X$
corresponding to the identity map $\{ {\rm Id}: H_{\X, T} \to H_{\X, T}\}\in$ $ \underline{H_{\X, T}}(H_{\X, T}):={\rm Mor}(H_{\X, T}, H_{\X, T})$. For any
$Y\in \underline{H_{\X, T}}(R)$ (so $Y$ is a closed subscheme in $\Sp(R\otimes_k S)$) we have $Y=\W_{\X,
T}\times_{H_{\X, T}} \Sp R$. In fact, the $k$-rational points of $\W$ are those pairs $(y, Y)$, where $Y\in
\underline{H_{\X, T}}(k)$ and $y\in \underline{Y}(k)$.

If $V$ is a finite dimensional $T$-module such that $k[V]^T=k$,
then $H^h_{V, T}$ is projective (see \cite[Corollary 1.2]{HS}) The
following lemma generalizes this statement.

\begin{lemma} \label{l0} Assume that $h(0)=1$.  Then the morphism $$p:H^h_{\X,T} \longrightarrow
\X/\!/T:=\Sp k[\X]^T$$ which assigns to any element $I\in
\underline{H^h_{\X,T}}(R)$ the morphism $k[\X]^T\to (R\otimes
k[\X]^T)/I^T\simeq R$, is projective.
\end{lemma}

\begin{proof}
Since we know that $H^h_{\X,T}$ quasiprojective, it is sufficient
to check that the valuative criterion of properness for $p$ is
satisfied. Let $S$ be the spectrum of a discrete valuation ring
$R$ with generic point $\eta$ and closed point $s$. We have to
show that any morphism $\phi_\eta: \eta\to H^h_{\X,T}$ such that
the composition $p\circ\phi:\eta\to \X/\!/T$ extends to a morphism
$S\to \X/\!/T$, extends to a unique morphism $S\to H^h_{\X,T}$.
Consider $Y_\eta = \eta\times_{H^h_{\X,T}} U^h_{\X,T}\subset
\eta\times \X$. By \cite[Prop. 9.7]{Har}, a closed subscheme
$Y\subset S\times \X$ such that $Y\times_S\eta=Y_\eta$, is flat
over $S$ if and only if $Y$ is the closure of $Y_\eta$ in $S\times
\X$. It follows that the desired extension  $Y\in
\underline{H^h_{\X,T}}(S)$ is unique. For the existence, we
consider  $Y:=\overline{Y_\eta}\subset S\times \X$. It remains to
show that that the fiber $Y_s$ is non-empty (then, by flatness, it
has the Hilbert function $h$). Indeed, we have the following
commutative diagram:

\begin{center}
\begin{picture}(120, 40)

\put(23,27){\vector(0,-1){10}} \put(82,27){\vector(0,-1){10}}

\put(5,35){$Y=\overline{Y_{\eta}} \subset\ \  S\times \X$}
\put(10,0){$S=\overline{\eta}\ \subset \ S\times \X/\!/T,$}


\end{picture}
\end{center}
where the morphism $Y\to S$ is the quotient by $T$, so it is
surjective.
\end{proof}

We prove the following lemma to treat one particular case of the
Hilbert scheme which we shall need later (see the corollary
below).

\begin{lemma} \label{l1}
Let $P$ be an $\N$-graded algebra: $P=\bigoplus_{r\geq 0} P_r$,
and

\medskip
\begin{center} $(*)$ \ \  there exists $r_0$ such that
$P_{r+1}=P_{1}P_{r}$ for any $r\geq r_0.$
\end{center}
\medskip

\noindent Consider the Hilbert scheme $H_{P}$ of the graded
algebra $P$ for the Hilbert function
$$h(r):= \left\{ \begin{array}{rl}
 1 & \mbox{if}\ \ r\geq 0, \\
 0 & \mbox{otherwise.}
\end{array}\right.$$ Let $R$ be an algebra and $Y=\Sp(R\otimes_k P/I)\in
\underline{H_{P}}(R)$. Then the projection $Y\to \Sp R$ is a
locally trivial bundle with fiber $\A^1$.
\end{lemma}

\begin{proof}
(1) Consider  the open subscheme in $\Sp P$ that is the complement
to the subscheme defined by the ideal $\bigoplus_{r>0} P_r$:
$$(\Sp P)_0=\{p\in \Sp P \ ; \
p\not\supseteq(\bigoplus_{r>0}P_r)\}$$ and the natural morphism
$$\psi: (\Sp P)_0\to \PP P.$$ Locally $\psi$ is given by the
embeddings of  algebras $(P_f)_0\subset P_f$, where $f\in P$ is
homogeneous, $\deg f> 0$ (it is clear that the corresponding
morphisms of affine schemes satisfy the compatibility conditions).
Note that $\psi$ is a locally trivial bundle with fiber $\G$.
Indeed, condition (*) implies that $\PP P$ is covered by open
affine subschemes $\Sp (P_{h})_0$, where $h\in P_1$, and for any
$h\in P_1$ we have $P_{h}=(P_{h})_0[h, h^{-1}]$.

(2) Consider $Y_0=Y\cap (\Sp R\times (\Sp P)_0)$. We have the
morphisms
$$Y_0\stackrel{\rho}{\longrightarrow} \PP (R\otimes_k P/I)\stackrel{\delta}{\longrightarrow} \Sp R.$$ Since $R\otimes_k P/I$
satisfies condition (*), by (1), it follows that $\rho$ is a
locally trivial bundle with fiber $\G$. So
 $\delta$ is an isomorphism. Consider the following morphism from $\Sp R$ to $\PP P$: $$\Sp R\cong \PP
(R\otimes_k P/I) \subset \Sp R\times \PP P
\stackrel{p}{\longrightarrow} \PP P,$$ where $p$ is the
projection.

(a) Note that $Y_0=(\Sp P)_0\times_{\PP P} \Sp R$. Indeed, locally
we have $$R\otimes_k P_f/I_f\simeq P_f\otimes_{(P_f)_0}
(R\otimes_k (P_f)_0/(I_f)_0),$$ where $f\in P$ is homogeneous of
positive degree.

(b) Consider $Y'=Y_0\times_{\G}\A^1$. Here $Y_0\times_{\G}\A^1$
denotes the categorical quotient $(Y_0\times \A^1)/\!/\G$, where
$\G$ acts on $\A^1$ as follows: $t\c s=t^{-1}s$, $t\in \G, s\in
\A^1$. Then  $Y'$ is a locally trivial bundle over $\Sp R$ with
fiber $\A^1$ and we have
 the natural morphism
$\eta : Y'\to Y$, which is locally given by the homomorphisms
$$R\otimes_k P/I\to \bigoplus_{r\geq 0} (R\otimes_k P_f/I_f)_r,$$
where $f\in P$ is homogeneous of positive degree. So we have a
commutative diagram:
\begin{center}
\begin{picture}(120, 40)
\put(25,30){$Y'$} \put(45,30){$\stackrel{\eta}{\longrightarrow}$}
\put(80,30){$Y$} \put(30,25){\vector(1,-1){15}}
\put(80,25){\vector(-1,-1){15}} \put(10,15){$$} \put(40,0){$\Sp
R.$} \put(22,13){$\alpha'$} \put(77,13){$\alpha$}
\end{picture}
\end{center}

Note that for any $r\geq 0,$ the corresponding homomorphism
$\alpha_*(\O_Y)_r\to \alpha'_*(\O_Y')_r$ is a surjective
homomorphism of locally free sheaves of $R$-modules of rank 1 and,
consequently, is an isomorphism. Thus $\eta$ is an isomorphism.
\end{proof}

The statement of the following corollary was given in
\cite[Section 5]{HS} with a proof for algebras generated by
elements of degree 1.

\begin{corollary} \label{cor1} With the notation of the previous lemma, the Hilbert scheme
$H_{P}$ is isomorphic to  $\PP P$.
\end{corollary}

\begin{proof}
We shall show that $\PP P$ represents the Hilbert functor
$\underline{H_{P}}$. For this we  prove that the tautological
bundle over $\PP P$ is the universal family, i. e., we are going
to prove the universal property for $E:=(\Sp P)_0\times_{\G}\A^1.$
Let $Y=\Sp(R\otimes_k P/I)\in \underline{H_{P}}(R)$. We have to
show that $Y=E\times_{\PP P} \Sp R$. Indeed, we have
$Y=Y_0\times_{\G}\A^1=((\Sp P)_0\times_{\PP P} \Sp
R)\times_{\G}\A^1=E\times_{\PP P} \Sp R.$
\end{proof}

Let us return to the case of an affine toric $\T$-variety  $\X$. We have
$$S=k[\X]=\bigoplus_{\nu\in\, \Omega}S_{\nu},$$ where
$\Omega \subset \Chi(\T)$ is a finitely generated monoid and $S_\nu$ is the subspace of $\T$-semiinvariant
functions of weight $\nu$ ($\dim S_{\nu} = 1$). Let $T\subset \T$ be a subtorus. We have a surjective
linear map $\pi : \Chi(\T)\to \Chi(T)$ given by the restriction.
 The action of $T$ on $\X$ arising from the action of $\T$ gives a grading  $$S=\bigoplus_{\chi\in \Sigma}S_{\chi},$$
 where $\Sigma=\pi(\Omega)$. We shall consider the following Hilbert function:
$$h(\chi):= \left\{ \begin{array}{rl}
 1 & \mbox{if}\ \ \chi\in \Sigma, \\
 0 & \mbox{otherwise.}
\end{array}\right. $$
Let $H_{\X, T}$ be the corresponding Hilbert scheme (we shall also denote it by $H_{S, T}$).
Note that all the ideals $I\in \underline{H_{V, T}}(k)$ are binomial (see \cite[Proposition 1.11]{ES}). If $x\in \X$ lies in the open $\T$-orbit, then we have the point $X:=\overline{T\c x}\in
\underline{H_{\X, T}}(k)$.

The group $\underline{\T}(R)$ acts on $\underline{H_{\X, T}}(R)$
 in the natural way. Namely, we have an action of $\underline{\T}(R)$ on $R\otimes_k S$ : for
 $f\in R\otimes_k S_{\nu}$, where $\nu\in \Omega$, and $t\in \underline{\T}(R) $
 let $t\c f= \nu(t)f$.
Hence for $I\in \underline{H_{\X, T}}(R)$ let $t\c I= \{t\c f \ ; \ f\in I\}$. These actions commute with  base
extensions, thus we have
  an action of $\T$ on $H_{\X, T}$. Since $T$ acts trivially, this yields an action of the torus $\T/T$.
   The universal family $\W_{\X, T}$ is invariant under the
diagonal action of $\T$
   on $H_{\X, T}\times \X$.

Let $H_0$ be the toric orbit closure $\overline{\T\c X}\subset H_{\X, T}$, and denote by $\W_0$ its preimage under
the projection $$p : \W_{\X, T}\to H_{\X, T}$$ (we consider $H_0$ and $\W_0$ with their structure of reduced
schemes).

\begin{proposition} \label{prop1} $(1)$ The stabilizer of $X$ under the action of $\T$ on $H_0$ is $T$. Moreover, $H_0$
 is a toric variety under the torus $\T/T$.

$(2)$ The orbit $\T\c X$ is open in $H_{\X, T}$. Consequently, $H_0$
is an irreducible component of $H_{\X, T}$.

$(3)$ $\W_0$ is a toric variety under the torus $\T$  $($and, consequently, $\W_0$ is an irreducible component of
$\W_{\X, T}$$)$. \end{proposition}

\begin{proof} (1) If $t\c X=X$ for $t\in \T$, then $t\c x\in T\c x$ and $t\in
T$. So we have only to show that $H_{\X, T}$ admits an open
covering by affine $T$-invariant charts. Indeed, let $\chi\in
\Sigma$. Then for any $I\in \underline{H_{\X, T}}(R)$ the locally
trivial $R$-module $(R\otimes k[\X]_\chi)/I_\chi$ defines a
morphism from $\Sp R$ to the projectivisation
$\Pr(k[\X]_{\chi}^*)=\PP ({\rm Sym}(k[\X]_{\chi}))$, where ${\rm
Sym}(k[\X]_{\chi})$ denotes the symmetric algebra. These maps
commute with base changes and, consequently, define a morphism
$p_\chi:H_{\X, T}\to \Pr(k[\X]_{\chi}^*)$. Note that $p_\chi$ is
$\T$-equivariant (the action of $\T$ on $\Pr(k[\X]_{\chi}^*)$ is
induced by the linear action of $\T$ on $k[\X]_{\chi}^*$). By
\cite[Proposition 3.2, Corollary 3.4]{HS}, it follows that there
exists a finite set of characters $\chi_1,\ldots, \chi_r\in
\Sigma$ such that the morphism
$$p\times p_{\chi_1}\times\ldots\times p_{\chi_r}: H_{\X, T}\to \X/\!/T\times\Pr(k[\X]_{\chi_1}^*)\times\ldots\times
\Pr(k[\X]_{\chi_r}^*)$$ is injective. Since the morphism $p$ is
projective (Lemma \ref{l0}), it follows that $p\times
p_{\chi_1}\times\ldots\times p_{\chi_r}$ is a closed embedding.
Since any $\Pr(k[\X]_{\chi_i}^*)$ admits an open covering by
$\T$-invariant affine charts, it follows that $H_{\X, T}$ does.

(2) We shall prove that $\T\c X$ is open in $H_{\X, T}$. Since the stabilizer of $X$ in $\T$ is $T$,
 it suffices to prove that $\dim T_X H_{\X, T} \leq \dim \T\c X=\dim \T - \dim T$, where $T_X H_{\X, T}$
 denotes the tangent space to $H_{\X, T}$ at $X$. By \cite[Prop. 1.6]{HS}, we have
 $$T_X H_{\X, T}=\Hom_{k[\X]}(I_X, k[X])_0.$$ This vector space is isomorphic to  $$\Hom_{k[\T]}(I_T, k[T])_0=
 \Hom_{k[\T]}(I_T/I_T^2, k[T])_0,$$
  where $I_T$ is the ideal of functions in $k[\T]$ vanishing on $T$. Indeed,
since $(I_X)_\chi\subset k[\X]_\chi(I_T)_0$, for any $\phi\in \Hom_{k[\T]}(I_T, k[T])_0$ we have $\phi(I_X)\subset k[\X]\phi((I_T)_0)=k[\X].$ Conversely, $I_T=k[\T]I_X$, so any $\phi\in \Hom_{k[\X]}(I_X, k[X])_0$ can be extended to a homomorphism of $k[\T]$-modules from $I_T$ to $k[T]$.

Further, we can choose  coordinates on $\T$ such that $$k[\T]=k[t_1,t_1^{-1},\dots,t_m,t_m^{-1},s_1,s_1^{-1},
...,s_r,s_r^{-1}],$$ where $r=\dim \T-\dim T$, and the ideal $I_T$ is generated by  $s_i-1$ for $i=1,\ldots, r$.
 The linear space $I_T$ is spanned by the elements $t_1^{a_1}\ldots t_n^{a_n}s_1^{b_1}\ldots s_m^{b_m}(s_i-1)$,
 where $a_i, b_j\in\Z$, and the projections of the elements $t_1^{a_1}\ldots t_n^{a_n}(s_i-1)$
 span the linear space $I_T/I_T^2$ (since $s_i(s_j-1)=(s_j-1)+(s_i-1)(s_j-1)$ and
 $s_i^{-1}(s_j-1)=(s_j-1)-s_i^{-1}(s_i-1)(s_j-1)$).
 Hence a homomorphism of $k[\T]$-modules from $I_T$ to $k[T]$ is uniquely determined by the images of
 $s_i-1$. Thus the dimension of the vector space of such homomorphisms of degree zero is not greater
 than $r$.

(3) Consider the restriction $p_0$ of $p$ to $\W_0$: $$p_0:
\W_0\to H_0.$$ This is a flat morphism. By Lemma \ref{lem} below
and \cite[Corollary 9.6]{Har}, the dimension of any irreducible
component $Z$ of $\W_0$ is equal to $\dim \T$. This implies that
$p_0(Z)=H_0$ and $Z\subset \overline{p^{-1}(\T\c X)}$. Thus
$\W_0=\overline{p^{-1}(\T\c X)}=\overline{\T\c(x, X)}$ is
irreducible and $\T\c(x, X)\subset \W_0$ is dense and,
consequently, open. Since $\W_0$ is a closed subscheme in
$H_0\times \X$, it follows that $\W_0$ admits an open covering by
affine $\T$-invariant charts.
\end{proof}

\begin{lemma} \label{lem} For any point $Y\in H_{\X, T}$, the dimension of any irreducible component of its
fiber $p^{-1}(Y)$ equals $\dim T$. \end{lemma}

\begin{proof} We denote  by $k(Y)$ the residue field of $Y\in H_{\X, T}$. Then we have $$p^{-1}(Y)=\Sp
k(Y)\times_{H_{\X, T}}\W_{\X, T}=\Sp L,$$ where $L$ is a coherent sheaf of $\Sigma$-graded $k(Y)$-algebras:
$$L= \bigoplus_{\chi\in\, \Sigma} L_\chi,$$ and $L_\chi:=k(Y)\otimes_{\O_{H_{\X, T}}}(\O_{\W_{\X, T}})_\chi$ is
isomorphic to $k(Y)$.  Let $$\Sigma_{red}:=\{\chi\in \Sigma \ ; \
L_\chi\ {\rm is\ not\ nilpotent}\}.$$ Note that
$\cone(\Sigma_{red})=\cone(\Sigma)$. Every point $Y\in H_{\X, T}$
gives us a subdivision of $\cone(\Sigma)$ into subcones, namely
two points $\chi, \chi' \in \Sigma_{red}$ lie in the same cone if
and only if $L_{\chi}L_{\chi'}\ne 0$. The irreducible components
$Z$ of $p^{-1}(Y)$  correspond to the maximal cones $C$ of this
subdivision: $$Z=\Sp( \bigoplus_{\chi\in\, \Sigma_{red}\cap
C}L_\chi).$$ Note that $\Sigma_{red}\cap C$ is a monoid. It
suffices to prove that the dimension of $Z$ is equal to $\dim T$.
We can extend the action of $T$ on $Z$ to an action of the torus
$T\times \Sp k(Y)$ (over the field $k(Y)$). Thus $Z$ is a toric
variety under the torus $T\times \Sp k(Y)$  and $\dim Z=\dim
C=\dim T$.
\end{proof}

\section{Fan of a toric Hilbert scheme}

Our aim is to describe the fans  of the toric varieties $H_0$ and $\W_0$.

Let us fix the notations. Recall that  $\X$ is an affine toric
variety under an action of a torus $\T$:  $$\X=\overline{\T\c
x_0},$$  $T\subset \T$ is a subtorus, and $$\pi : \Chi(\T)\to
\Chi(T)$$ is the restriction map. Fix isomorphisms $\T\simeq \G^n,
T\simeq \G^r$, this gives us  a basis in $k[\T]$ and $k[T]$:
$$k[\T]=\bigoplus_{\nu\in\Chi(\T)}kt^\nu,\ \ \ \ \ \ \ \ \ \ \
k[T]=\bigoplus_{\chi\in\Chi(T)}kt^\chi.$$
 We
denote by $X$ the $T$-orbit closure $\overline{T\c x_0}$  and
$I_X\subset k[\X]$ denotes the corresponding ideal. Also, we have
$$\! \! \! S=k[\X]=\bigoplus_{\nu\in \Omega}kt^\nu,\ \ \ \ \ \ \ \
k[X]=\bigoplus_{\chi\in \Sigma}kt^\chi.$$ The restriction
homomorphisme $k[\X]\to k[X]$ is given by $t^\nu\to t^{\pi(\nu)}$
and its kernel $I_X$ is generated by all
 the binomials of the form $t^{\nu_1}-t^{\nu_2}$ such that $\pi(\nu_1)=\pi(\nu_2)$ (see \cite[Lemma 4.1]{St}).

Let us recall some definitions concerning convex polyhedra. They are taken from \cite{St}, which we shall use as a general reference on convex polyhedra.

\begin{definition} Let $P$ be a convex polyhedron in a vector space $V$. For any face $F$ of $P$ the {\it normal cone} $N_F(P)$ is
the following cone in the dual vector space $V^*$: $$N_F(P):=\{l\in V^*\ ; \ l(v-v')\geq 0 {\rm \ for\  all\ } v\in P, v'\in F \}.$$
 The {\it normal fan} $N(P)$ of $P$ is the fan whose cones are normal cones to the faces of $P$.
\end{definition}

\begin{definition}
The {\it recession cone} of a polyhedron $P\subset V$ is
the set of those vectors $v\in V$ such that $u+v\in P$ for any $u\in P$.
\end{definition}

\begin{definition} A fan $\C_1$ is a {\it refinement} of a fan $\C_2$ if any cone of $\C_1$ is
contained in some cone of $\C_2$.
\end{definition}

\begin{definition} We say that two polyhedra $P_1, P_2\subset \Chi(\T)_{\R}$ are {\it equivalent}
if they have the same normal fan.
\end{definition}

We fix an open $\T$-equivariant embedding of $\T/T$ (resp. of
$\T$) in $H_0$ (resp. in $W_0$) such that the image of $eT$ (resp.
of $e$) is $X$ (resp. $(X, x_0)$), where $e$ is the unit in $\T$.

\begin{theorem} \label{th1}  $(1)$ The fan $\C_{H_0}\subset \Lambda(\T)_{\R}$ of the toric $\T/T$-variety $H_0$
is the coarsest common refinement of the normal fans $\C_{\chi}$ of the polyhedra
$$P_{\chi}:={\rm conv}(\pi^{-1}(\chi)\cap \Omega)\subset\Chi(\T)_{\R},$$ where $\chi\in \Sigma$.

$(2)$ The fan of $\W_0$ is the coarsest common refinement of the
fans $\C_{H_0}$ and $N(\cone(\Omega))$.
\end{theorem}

\begin{remark} We can consider fans in $\Lambda(\T/T)_{\R}$ as fans in $\Lambda(\T)_{\R}$ whose cones
contain $\Lambda(T)_{\R}$. In particular, we view the fan of the toric $\T/T$-variety $H_0$ as a fan in
$\Lambda(\T)_{\R}$. \end{remark}

\begin{proof} Let us recall that a one-parameter subgroup $\lambda\in\Lambda(\T)_{\R}$ belongs to the
support of the fan $C_{H_0}$ if and only if there exists a limit
of $X\in H_0$ under $\lambda$. Further, one-parameter subgroups
$\lambda, \lambda'\in \Lambda(\T)_{\R}$ lie in the interior of the
same cone of $C_{H_0}$ if and only if they define the same limit
of $X\in H_0$.

We shall calculate the limit of $X$ under a one-parameter subgroup
$\lambda\in \Lambda(\T)$. Consider the closed embedding $$\G\times
X\subset \G\times \X,$$
$$(s,x)\to (s,\lambda(s)\c x).$$ Let $\Xi$ be the closure of the
image of this embedding in $\A^1\times \X$ (so $\Xi$ is a
variety). Since the projection $p_{\A^1}:\Xi\to \A^1$ is a flat
morphism, we have a morphism $\A^1\to H_{\X, T}$ such that
$\Xi=\W_{\X, T}\times_{H_{\X, T}}\A^1$. Thus the limit $X_\lambda$
of $X$ under $\l$ is equal to the fiber of $p_{\A^1}$ over 0 if
this fiber is non-empty and the limit does not exist otherwise.
Consider the commutative diagram:
$$
\begin{array}{ccc}
\Xi & \supset & \G \times X \\
\cap && \cap \\
\A^1\times \X & \supset & \G \times \X. \\
\end{array}
$$
We have the corresponding homomorphisms of algebras :
$$
\begin{array}{ccc}
k[\Xi] & \hookrightarrow & k[\G \times X] \\
\uparrow && \uparrow \\
k[\A^1 \times \X] & \hookrightarrow & k[\G \times \X], \\
\end{array}
$$
where the vertical maps are surjective.

Denote by $s$ the coordinate in $\A^1$. Then the homomorphism
$k[\G \times \X]\twoheadrightarrow k[\G \times X]$  is given by
$s\to s$ and $t^\nu\to s^{\langle\lambda,
\nu\rangle}t^{\pi(\nu)}$. Thus, the vector subspace $k[\Xi]\subset
k[\G \times X]$ is generated by the elements of the form
$s^mt^{\pi(\nu)}$, where $m\geq \langle\lambda, \nu\rangle$,
$\nu\in\Omega$. The fiber $p_{\A^1}^{-1}(0)$ is empty if and only
if the ideal $sk[\Xi]$ contains 1. Thus $\lambda$ belongs to the
support of the fan $\C_{H_0}$ if and only if
$\langle\lambda,\nu\rangle\leq 0$ for any $\nu\in \pi^{-1}(0)\cap
\Omega$. Since any $k[\X]_\chi$ is a finitely generated
$k[\X]_0$-module, this is equivalent to say that $\lambda$ attains
its minimum on $\pi^{-1}(\chi)\cap \Omega$ for any $\chi\in
\Sigma$. In this case
$$k[X_\lambda]=\bigoplus_{\chi\in\, \Sigma}\, ks^{n_{\l}(\chi)}t^{\chi},$$ where
$$n_{\l}(\chi):=\min_{\nu\in \, \pi^{-1}(\chi)\cap \Omega}
\langle\lambda,\nu\rangle.$$  The product
$s^{n_{\l}(\chi_1)}t^{\chi_1}s^{n_{\l}(\chi_2)}t^{\chi_2}=s^{n_{\l}(\chi_1)+n_{\l}(\chi_2)}t^{\chi_1+\chi_2}$
equals zero if and only if $n_{\l}(\chi_1)+n_{\l}(\chi_2)>
n_{\l}(\chi_1+\chi_2)$. The embedding of $X_\lambda$ in $\X$ is
given by the homomorphism of algebras $k[\X]\to k[X_\lambda]$,
where $t^\nu$, $\nu\in \Omega,$ maps to
$s^{n_{\l}(\pi(\nu))}t^{\pi(\nu)}$ if $\langle\lambda,
\nu\rangle=n_{\l}(\chi)$, and to 0 otherwise. We denote by
$I_\lambda$ the kernel of this homomorphism. Hence we see that
one-parameter subgroups $\lambda_1$ and $\lambda_2$ define the
same limit if and only if $I_{\lambda_1}=I_{\lambda_2}$. This
holds if and only if $\lambda_1$ and $\lambda_2$ attain the
minimum over $\pi^{-1}(\chi)\cap \Omega$ at the same point for any
$\chi\in \Sigma$ or, equivalently,  $\l_1$ and $\l_2$ lie in the
interior of the same cone of $N(P_{\chi})$ for any
$\chi\in\Sigma$.

(2) Since $\W_0=\overline{\T\c(X, x)}\subset H_0\times \X$, the second statement is evident.

\end{proof}

The statement below follows directly from the description of the
limit of $X$ under $\lambda$ in the proof of the theorem.

\begin{remark} Let $\prec_\lambda$ be the preorder on $\Chi(\T)$ such that $\nu_1\prec_\lambda\nu_2$ if
 $\langle\lambda, \nu_1\rangle\leq \langle\lambda, \nu_2\rangle$. For any $f=\sum f_{\nu_i}, f_{\nu_i}\in k[\X]_{\nu_i}$,  denote by ${\rm in}_{\lambda}(f)$ the sum of $f_{\nu_i}\!$'s where $\nu_i$ is maximal with respect to $\prec_\lambda$. Then the limit  of  $I_X\in H_0$ under $\lambda$ exists if and only if $\langle\lambda, \nu\rangle\geq 0$ for any $\nu\in \pi^{-1}(0)\cap\Omega$. In this case the limit is the ideal ${\rm in}_{\lambda}(I_X)$ generated by all ${\rm  in}_{\lambda}(f), f\in I_X$.

\end{remark}

\begin{example} Let $\X=\A^n$, $\T=\G^n$ act on $\A^n$ by rescaling of coordinates, $T=\G$, and let the $\Chi(T)$-grading of $k[x_1,\ldots, x_n]$ be positive.

(1) Consider the case $n=3$. It was proved by Arnold, Korkina, Post, Roelfols (see, for example, \cite[Theorem
10.2]{St}), that any ideal $I\in \underline{H_{\A^n, T}}(k)$ is of the form $t\c in_{\lambda}(I_X)$ for some $t\in \G^n$ and
$\l\in \Lambda(\G^n)$. This means that in this case the toric Hilbert scheme is irreducible.

(2) Let $n=4$ and $\chi_1=1, \chi_2=3, \chi_3=4, \chi_4=7$. Then the toric Hilbert scheme is reducible.
Moreover, in $H_{\A^n, T}$ there are infinitely many orbits of $\G^n$ (see \cite[Theorem 10.4]{St}).
\end{example}

\begin{proposition}\label{cor}
$(1)$ The support of any $\C_\chi$ is the cone generated by those
one-parameter subgroups $\lambda$ that $\langle\lambda,
\nu\rangle\geq 0$ for any $\nu\in \pi^{-1}(0)\cap \Omega$. In
particular, the grading of $S$ by $\Chi(T)$ is positive if and
only if this support is the whole space $\Lambda(\T)_{\R}$, i.e.,
any polyhedron $P_{\chi}$ is a polytope. This holds if and only if
$H_0$ is projective.

$(2)$ There are only  finitely many non-equivalent polyhedra $P_{\chi}$ for $\chi\in \Sigma$. Hence $\C_{H_0}$ is the
normal fan of the Minkowski sum of representatives of the equivalence classes $($we denote this sum by
$P_{H_0})$.
\end{proposition}

\begin{proof}
(1) First note that $P_0$ is a cone and its normal cone $\C_0$ is generated by those one-parameter subgroups
$\lambda$ that $\langle\l, \nu\rangle\geq \langle\l, 0\rangle = 0$ for any $\nu\in \pi^{-1}(0)\cap
\Omega$. Further, note that the recession cone of any $P_\chi$ is $P_0$. Indeed, $S_\chi$ is a finitely generated
$S_0$-module. Let $\mu_1,\ldots,\mu_d\in\Chi(\T)$ be the weights of a set of $\T$-semiinvariant generators.
Then
$$P_\chi={\rm conv}(\bigcup_{i=1}^d(\mu_i+P_0))={\rm conv}(\mu_1,\ldots, \mu_d)+P_0.$$ It follows
that the support of $\C_\chi$ is $\C_0$.

If the support of $\C_{H_0}$ is not $\Lambda(\T)_{\R}$, then $H_0$ is not complete and, consequently, is not
projective. Conversely, if the grading is positive, then the Hilbert scheme $H_{\X, T}$ is projective, and $H_0$ is projective.

(2) There are  only finitely many  fans $\C$ such that $\C_{H_0}$ is a refinement of $\C$ and the supports of
$\C$ and $\C_{H_0}$ coincide.\end{proof}

\begin{remark} By \cite[Theorem 7.15]{St}, it follows that in the case when $\X=\A^n$, $\T=\G^n$ acts by rescaling of coordinates, and the $\Chi(T)$-grading of $k[\X]$ is positive,  the polytope
$P_{H_0}$ is equivalent to the Minkowski sum of $P_{\chi}$ corresponding to the weights $\chi$ of the elements of the universal
Gr\"{o}bner basis of $I_{\X}$.\end{remark}

Let $\X$ be normal. Now we are going to give a precise description of those characters $\chi\in \Sigma$ having
equivalent polyhedra $P_\chi$. Recall that we have a homomorphism of lattices $\pi : \Chi(\T)\to \Chi(T)$, a finitely
generated monoid $\Omega\subset \Chi(\T)$ such that $\Omega=\cone(\Omega)\cap \Chi(\T)$, and we put
$\Sigma=\pi(\Omega)$. To any point $\chi\in \Sigma$ we associate the polyhedron $$P_\chi ={\rm conv}(\pi^{-1}(\chi)\cap
\Omega)\subset \Chi(\T)_{\R}.$$ Two points $\chi, \chi'\in \Sigma$ are said to be {\it equivalent} if the corresponding polyhedra
$P_{\chi}$ and $P_{\chi'}$ are equivalent. The question is to describe equivalence classes constructively.

Denote by $\pi_{\R}$ the linear map induced by $\pi$ : $$\pi_{\R} : \Chi(\T)_{\R}\to \Chi(T)_{\R}.$$ Let
$\C^{\R}_\chi$ denote the normal fan to the polyhedron $$P^{\R}_\chi:=\pi_{\R}^{-1}(\chi)\cap\cone(\Omega).$$

\begin{definition} \label{cell} (see \cite{CM}) The {\it cell decomposition of} $\cone(\Sigma)$ {\it induced by} $\pi_{\R}$ is the subdivision of  $\cone(\Sigma)$ into the following set of cones : the characters $\chi$ and $\chi'$ lie in
the interior of the same cone of this decomposition if and only if the set of those faces of $\Omega_{\R_+}$ whose images under $\pi_{\R}$ contain $\chi$ coincides with the set of such faces for $\chi'$.
\end{definition}

\begin{remark}
Note that the cell decomposition of $\cone(\Sigma)$ induced by $\pi_{\R}$ coincides with the subdivision by GIT-cones (\cite[Section 2]{BH}).
\end{remark}

Note that if $\chi$ lies in the interior of a cone $\sigma$ of the cell decomposition and $\chi'\in
\sigma$, then $\C^{\R}_\chi$ refines $\C^{\R}_{\chi'}$. In particular, the polyhedra $P^{\R}_{\chi}$ corresponding to
interior points $\chi$ of $\sigma$ are equivalent. Let $P_{\R}$ denote the Minkowski sum of $P^{\R}_{\chi}$ for
representatives of interior points for all cones of the cell decomposition and let $\C_{\R}$ denote the normal fan to
$P_{\R}$ (note that in the Minkowski sum it  suffices to take representatives of interior points for the
maximal cones of the cell decomposition).

\begin{remark} In \cite{CM} the fan $\C_{\R}$ is called the {\it fiber fan} by analogy with the normal
fan of the fiber polytope for a linear projection of polytopes (see \cite{BS}).
\end{remark}

\begin{definition}\label{int} (See \cite[Definition 5.4]{HS}.) A character $\chi\in \Sigma$ is {\it integral} if the inclusion of
the convex polyhedra $P_\chi\subseteq P^{\R}_\chi$ is an equality.
\end{definition}

We shall denote by $\Sigma^{\rm int}_{\X}$ the set of integral characters. The following proposition gives us an algorithm for computing the fan of $H_0$.

\begin{proposition} \label{prop2} For any cone $\sigma$ of the cell decomposition of $\cone(\Sigma)$ induced by $\pi$
let $\mu_1,\ldots, \mu_r$ be generators of the monoid $\sigma\cap \Sigma$ and let $c_1,\ldots, c_r\in \N$ be
such that $c_i\mu_i$ are integral, $i=1,\ldots, r$. Then, the polyhedra $P_\chi$, where
$\chi=\sum_{i=1}^rd_i\mu_i$ and $0< d_i < l(\sigma)c_i$, form representatives of all equivalence classes of points
in $\sigma$ up to Minkowski sum with $P_{\R}$. Here $l(\sigma)$ is the number of vertices of $P^{\R}_\chi$ for $\chi$ lying in the interior of $\sigma$.

Hence $P_{H_0}$ is the Minkowski sum of such
representatives for all (maximal) cones $\sigma$ of the cell decomposition of $\cone(\Sigma)$ induced by $\pi_{\R}$.
\end{proposition}

\begin{proof}

Consider a point $\chi$ lying in the interior of $\sigma$ and the
corresponding polyhedron $P^{\R}_{\chi}$. For any vertex $v$ of
$P^{\R}_{\chi}$ there exists a unique minimal face $F$ of
$\cone(\Omega)$ such that $F\cap P^{\R}_{\chi}=\{v\}$ (indeed,
since $P^{\R}_{\chi}$ is the intersection of $\cone(\Omega)$ with
the affine subspace $\pi_{\R}^{-1}(\chi)$, it follows that any
face  of $P^{\R}_{\chi}$ is the intersection of
$\pi_{\R}^{-1}(\chi)$ with some face of $\cone(\Omega)$). Let
$v^\chi_1, \ldots, v^\chi_{l(\sigma)}\in \Chi(\T)_{\R}$ be the
vertices of $P^{\R}_{\chi}$ and let $F_1^\sigma, \ldots,
F_{l(\sigma)}^\sigma$ be the corresponding faces (the set of such
faces does not depend on a point $\chi$ in the interior of
$\sigma$). Note also that the intersection $F_i\cap P^{\R}_{\chi}$
is a vertex of $P^{\R}_{\chi}$ for any $\chi\in\sigma$. Just as
above, we denote this vertex by $v^\chi_i$. For two vectors $u,
u'\in \Chi(\T)_{\R}$ we say $u\prec u'$ if $u'-u\in
\cone(\Omega)$.

Let us show that if $\chi=\sum_{i=1}^rd_i\mu_i$ lies in the
interior of $\sigma$ and there exists $i$ such that $d_i\ge c_i
l(\sigma)$, then
 $$(*)\ \ \ \ \ \ \ P_{\chi}=P_{\chi-c_i\mu_i}+P_{c_i\mu_i }.$$
Indeed, the inclusion $P_{\chi-c_i\mu_i}+P_{c_i\mu_i}\subseteq
P_{\chi}$ is evident. For the converse, it is sufficient to show
that $P_{\chi}\cap \Omega\subset P_{\chi-c_i\mu_i}+P_{c_i\mu_i}$.
Note that $P_{\chi}\cap \Omega=P^\R_{\chi}\cap\Omega$. Denote by
$D_\chi$ the convex hull of the $v^\chi_i$,
$i=1,\ldots,l(\sigma)$. By Proposition \ref{cor} (1),
$P^\R_{\chi}=D_\chi+P_0$.
 Then for
any $v\in P_{\chi}\cap \Omega$ we have $v=u+v_0$ for some $v_0\in
P_0, u\in D_\chi$ and
$u=\sum_{j=1}^{l(\sigma)}q_jv^{\chi}_j$ for some $q_j\geq 0$ such
that $\sum_{j=1}^{l(\sigma)}q_j=1$. There exists $j$ such that
$q_j\geq 1/l(\sigma)$. Hence $v\succ
q_jv^{\chi}_j=q_j(v^{\chi-c_i\mu_i}_j+v^{c_i\mu_i}_j)\succ
v^{c_i\mu_i}_j$. Thus $v-v^{c_i\mu_i}_j\in \cone(\Omega)\cap
\Chi(\T)=\Omega$.

In particular, this implies that for any $\chi$ in the interior of $\sigma$ the polyhedron $P_{\chi}+P_{\R}$ is equivalent
to $P_{\chi'}+P_{\R}$ for some $\chi'=\sum_{i=1}^rd_i\mu_i$ such that $d_i< c_i l(\sigma)$ for any $i$.
The second statement of the proposition is evident.
\end{proof}

\begin{corollary} With the preceding notation,
if $\chi=\sum_{i=1}^rd_i\mu_i$ lies in the interior of $\sigma$ and there exists $i$ such that $d_i\ge c_i l(\sigma)$, then $P_{\chi}$ is equivalent to $P_{\chi+c_i\mu_i}$.
\end{corollary}

\begin{proof}
By $(*)$, it follows that $P_{\chi+c_i\mu_i}=P_{\chi-c_i\mu_i}+2P_{c_i\mu_i}$ is equivalent to $P_{\chi}$.

\end{proof}

\begin{example} Let $\X=\A^n$, $\T =\G^n$ act be rescaling of coordinates, and let $T=\G$ act on $\A^n$ with characters $\chi_1,\ldots,$ $\chi_n\in \Z$. Then
$\Omega\subset \Z^n$ is the set of vectors with integral non-positive coordinates, and
$\Sigma\subset \Z$ is the monoid generated by $-\chi_i$. Moreover, $\Sigma =(\Sigma \cap \Z_+) \cup (\Sigma \cap
\Z_{-})$ is the subdivision of $\Sigma$ induced by $\pi$. Let $n_+$ and $n_-$ be the numbers of positive and
negative $\chi_i$ respectively. A number $\chi\in \Z_+$ (resp. $\Z_-$) is integral (in the sense of Definition \ref{int}) if and only if $\chi$ is divisible
by any $\chi_i<0$ (resp. $>0$). Let $\chi_+$ (resp. $\chi_-$) be the least common (positive) multiple of all positive
(resp. negative) $\chi_i$. Then $P_{H_0}$ is the Minkowski sum of polyhedra $P_\chi$ for $-n_+\chi_+<\chi<n_-\chi_-$.
\end{example}

\section{Toric Chow morphism}

We are going to describe the toric Chow morphism from the Hilbert scheme to the inverse limit of GIT quotients
$\X/_{\!\!\chi} T$. In \cite[Section 5]{HS} the toric Chow
morphism was constructed in the case when $\X=\A^n$ is a $T$-module. We generalize this to the case of a normal
affine toric $\T$-variety $\X$.

In this section we fix a $\T$-equivariant closed embedding $\X\hookrightarrow V$, where $V$ is a finite-dimensional $\T$-module
such that $\X$ is not contained in a proper $\T$-submodule.
We use the notations of the previous sections. Let
$$S^{(\chi)}:=\bigoplus_{r=0}^{\infty} S_{r\chi},$$ and let $$\X/_{\!\!\chi} T:=\PP\ S^{(\chi)}$$ be the GIT quotient. In particular,
$\X/_{\! 0}T=\X/\!/T=\Sp (S_0)$. Notice also that $\X/_{\!\!\chi} T = \X^{ss}_{\chi}\!/\!/T$, where $$\X^{ss}_{\chi}:=\{x\in \X ; f(x)\ne 0  {\rm\ for\ some\
homogeneous\ }  f\in S^{(\chi)}\}.$$ If $\chi$ lies in the interior of $\cone(\Sigma)$, then $\X/_{\!\!\chi} T$ is a
normal toric $\T/T$-variety whose fan is $\C^{\R}_\chi$, the normal fan to the
polyhedron $P^{\R}_\chi$.

It is easy to see that for any $\chi_1, \chi_2\in \Sigma$, the
inclusion $\X^{ss}_{\chi_1}\subseteq \X^{ss}_{\chi_2}$ holds if
and only if $\chi_1$ belongs to the cone of the cell decomposition
of $\cone(\Sigma)$ induced by $\pi$ (see Definition~\ref{cell})
containing $\chi_2$ in its interior. We consider the morphisms
between GIT-quotients $\X/_{\!\!\chi} T$ induced by inclusions
between $\X^{ss}_{\chi}$, where $\chi\in \Sigma$. So the
GIT-quotients $\X/_{\!\!\chi} T$ form a finite inverse system with
$\X/\!/T$ sitting at the end. Consider the inverse limit
$$\X/_{\!\! C}T:=\lim_{\longleftarrow}\{\X/_{\!\!\chi} T ; \chi
{\rm\ lies\ in\ the\ interior\  of\ } \Sigma\}.$$ It is a closed
subscheme in the product $\X/_{\!\!\chi_1} T\times\ldots\times
\X/_{\!\!\chi_r} T$, where $\chi_1,\ldots,\chi_r$ are
representatives of interior points of all maximal cones of the
cell decomposition of $\cone(\Sigma)$ induced by $\pi$. Note also
that $\X/_{\!\! C}T$ is a closed subscheme in $V/_{\!\! C}T$.

\begin{definition}
The {\it main component} $(\X/_{\!\! C}T)_0$ of the inverse limit  $\X/_{\!\! C}T$ is  the closure of the image of the map $\T\to \X/_{\!\! C}T$ induced by the maps $\T\to \X/_{\!\!\chi} T$, where $\chi$ lies in the interior of $\Sigma$.
\end{definition}

By \cite[Proposition 3.8]{CM}, it follows that the main component $(\X/_{\!\! C}T)_0$ is an irreducible component of  $\X/_{\!\! C}T$ which satisfies the following universal property : given a $\T/T$-variety $Y$ containing an irreducible component $Y_0$ such that $Y_0$ is a toric $\T/T$-variety, and given $\T/T$-equivariant morphisms $\phi_\chi: Y\to \X/_{\!\!\chi} T$, where $\chi$ lies in the interior of $\Sigma$, such that the $\phi_\chi$ induce birational morphisms $Y_0\to \X/_{\!\!\chi} T$  and the $\phi_\chi$ are compatible with the morphisms of the inverse system (so the $\phi_\chi$ give a morphism $\phi: Y\to \X/_{\!\! C}T$); then restricting the morphism $\phi$ to $Y_0$ we have a birational morphism of toric $\T/T$-varieties $Y_0\to (\X/_{\!\! C}T)_0$.

\begin{remark} By \cite[Proposition 3.10]{CM}, it follows that the fan of  $(\X/_{\!\! C}T)_0$ is $\C_\R$, the
maximal common refinement of all the normal fans to the polyhedra $P^{\R}_\chi$, $\chi\in \Sigma$. Since every
character $\chi\in \Sigma$ has some integral positive multiple $c\chi\in \Sigma^{\rm int}_\X$ ($c\in \N$), the fan $\C_{H_0}$ is a refinement
of the fan $\C_\R$.
\end{remark}

The following example shows that $\C_{H_0}$ and $\C_\R$ do not always coincide.

\begin{example} \label{ex} Let $\X=\A^3$, $\T=\G^3$ act by rescaling of coordinates, and let $T=\G$ act by
$t(x_1,x_2,x_3)=(tx_1,tx_2,t^2x_3)$.

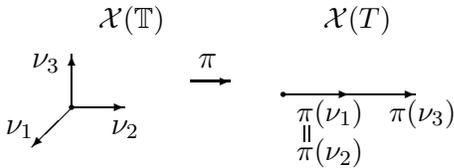
\begin{picture}(250, 70)
{\small \put(35, 50){$\Chi(\T)$}}

\put(25, 20){\circle*{2}}

\put(25, 20){\vector(-1, -1){15}} \put(0, 10){$\nu_1$}

\put(25, 20){\vector(1, 0){20}} \put(40, 10){$\nu_2$}

\put(10, 35){$\nu_3$} \put(25, 20){\vector(0, 1){20}}

\put(70, 30){\vector(1, 0){15}} \put(73, 35){$\pi$}







{\small \put(120, 50){$\Chi(T)$}

\put(110, 15){$\pi(\nu_1)$} \put(110, 0){$\pi(\nu_2)$} \put(113, 7){$\line(0,1){6}$} \put(115,
7){$\line(0,1){6}$}

\put(145, 15){$\pi(\nu_3)$}}




\put(105, 25){\vector(1, 0){25}}

\put(105, 25){\vector(1, 0){50}}

\put(105, 25){\circle*{2}}

\end{picture}
\medskip

The Hilbert scheme $H_{\A^3, T}$ is the closed subscheme in $\Pr^1\times\Pr^3$ defined by the equations
$z_1w_3-z_2w_1=0$ and $z_1w_2-z_2w_3=0$ (where $z_1,z_2$ and $w_1,w_2,w_3,w_4$ are homogeneous coordinates in $\Pr^1$ and $\Pr^3$
respectively). The integral (in the sense of Definition \ref{int}) degrees are even.
The fan $\C_{H_0}$ consists of the following cones:
$$\!\!\!\! \! \!\! \!\! \R_+(e_1+e_2)+\R_+e_2,$$
$$ \R_+(e_1+e_2)+\R_+(-e_2),$$
$$\!\!\!\! \! \! \!\!\! \R_+(e_2-e_1)+\R_+e_2,$$
$$ \R_+(e_2-e_1)+\R_+(-e_2),$$
where $e_1=\nu_1^*+\nu_3^*, e_2=-\nu_3^*$ is a basis of
$\Lambda(\T/T)$. The inverse limit of GIT-quotients is $\A^3/_{\!\! C}T=\PP\ k[x_1,x_2,x_3]$ (where $k[x_1,x_2,x_3]$ is graded by the
weights of $T$), and its fan $\C_\R$ consists of the following cones: $$\! \! \! \R_+(e_1+e_2)+\R_+(-e_2),$$
$$\! \! \! \R_+(e_2-e_1)+\R_+(-e_2),$$ $$\ \  \R_+(e_1+e_2)+\R_+(e_2-e_1).$$

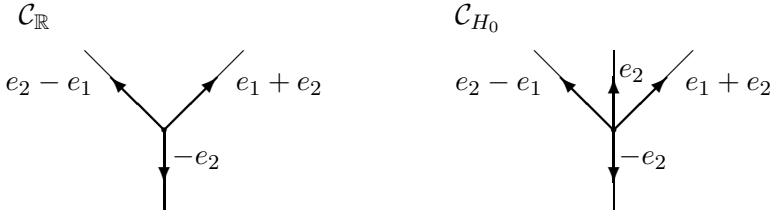
\begin{picture}(250, 90)
{\small \put(5, 80){$\C_\R$}

\put(60, 40){\circle*{2}}

\thinlines

\put(60, 40){\line(0, -1){30}} \put(60, 40){\line(1, 1){30}} \put(60, 40){\line(-1, 1){30}}

\thicklines

\put(60, 40){\vector(0, -1){20}} \put(60, 40){\vector(1, 1){20}} \put(60, 40){\vector(-1, 1){20}}

\put(63, 27){$-e_2$} \put(0, 55){$e_2-e_1$} \put(87, 55){$e_1+e_2$}

}

{\small \put(170, 80){$\C_{H_0}$}
\put(230, 40){\circle*{2}}

\thinlines

\put(230, 40){\line(0, -1){30}} \put(230, 40){\line(1, 1){30}} \put(230, 40){\line(-1, 1){30}}
\put(230, 40){\line(0, 1){30}}

\thicklines

\put(230, 40){\vector(0, -1){20}} \put(230, 40){\vector(1, 1){20}} \put(230, 40){\vector(-1, 1){20}}
\put(230, 40){\vector(0, 1){20}}

\put(232, 27){$-e_2$} \put(170, 55){$e_2-e_1$} \put(257, 55){$e_1+e_2$}
\put(232, 60){$e_2$}

}

\end{picture}

\end{example}

By the statement $(*)$ from the proof of Proposition \ref{prop2}, it follows that if a character $\chi\in \Sigma$ is integral, then there exists $r_0$ such that
 $S_{(r+1)\chi}=S_{\chi}S_{r\chi}$ for all $r\geq r_0$. Thus
Corollary \ref{cor1} implies that $$H_{S^{(\chi)},T}=\PP
S^{(\chi)}=\X/_{\!\!\chi} T,$$ for any $\chi\in \Sigma^{\rm
int}_{\X}$.

For any subset $D\subset \Sigma$ we can consider the restriction of the Hilbert scheme $H_{\X, T}$ on degrees $D$, that is, the quasiprojective scheme $H^D_{\X, T}$ representing the covariant functor $$\underline{H^D_{\X, T}} : \underline{k-Alg}\to \underline{Set}$$ such that
$\underline{H^D_{\X, T}}(R)$ is the set of families  $\{L_\chi\}_{\chi\in D}$, where $L_\chi\subset R\otimes_k S_\chi$ is an $R$-submodule, such that
$ (R\otimes_k S_\chi)/L_\chi$ is a locally free $R$-module  of rank 1 and
$fL_{\chi_2}\subset L_{\chi_1}$ for  any  $\chi_1, \chi_2 \in D$ and  any $f\in S_{\chi_1-\chi_2}$
(see \cite[Section 2]{HS}). In particular, $H^\Sigma_{\X, T}=H_{\X, T}$ and  $H_{S^{(\chi)}, T}=H^{D^\chi}_{\X, T},$ where $D^\chi:=\{c\chi\ ; \ c\in \Z_+\}$. Note also that $H^D_{\X, T}$ is a closed subscheme of $H^D_{V, T}$.  For any $D\subset \Sigma$ we have a degree restriction morphism $H_{\X, T}\to H^D_{\X, T}$. In particular, we have canonical morphisms $$\phi_{\X}^\chi:H_{\X, T}\to \X/_{\!\!\chi} T.$$

The following theorem was proved in \cite[Theorem 5.6]{HS} for the case when $\X$ is a finite-dimensional $T$-module.

\begin{theorem} \label{chow} Let $ H^{\rm int}_{\X, T}:= H^{\Sigma^{\rm int}_\X}_{\X, T}$ \ be the toric Hilbert scheme restricted to the set of integral degrees. Then there is a canonical morphism
$$\phi_{\X}^{\rm int} : H^{\rm int}_{\X, T}\to \X/_{\!\! C} T$$ which induces an isomorphism of the corresponding reduced schemes. In particular, composing $\phi_{\X}^{\rm int}$ with the degree restriction morphism, we obtain a canonical Chow morphism from the toric Hilbert scheme to the inverse limit of the GIT quotients  $$\phi_{\X} : H_{\X, T}\to \X/_{\!\! C} T.$$\end{theorem}

\begin{proof}
As in \cite[Lemma 5.7]{HS}, we see that the morphisms
$\phi_\X^\chi$ satisfy the compatibility  conditions for $\chi\in
\Sigma^{\rm int}_\X$ and, consequently, give a canonical morphism
$$ H_{\X, T} \to H^{\rm int}_{\X, T} \stackrel{\phi_{\X}^{\rm
int}}{\longrightarrow} \X/_{\!\! C} T.$$ Further, note that for
any algebra $R$ the morphism  $$\underline{\phi_{\X}^{\rm
int}}(R):\underline{H^{\rm int}_{\X, T}}(R)\to \underline{
\X/\!\!_C} T(R)$$ is injective (since
$H_{S^{(\chi)},T}=\X/_{\!\!\chi} T,$ we view any element of
$\underline{ \X/\!\!_C} T(R)$ as a family of $R$-submodules
$\{I_\chi\subset R\otimes S_\chi\}_{\chi\in \Sigma^{\rm int}_\X}$
such that $(R\otimes S_\chi)/I_\chi$ is a locally free $R$-module
of rank 1, $I^{(\chi)}:=\bigoplus_{n\geq 0} I_{n\chi}$ is an ideal
in $R\otimes S^{(\chi)}$, so it defines a point of $
\underline{H_{S^{(\chi)},T}}(R)_{\chi\in \Sigma^{\rm int}_\X}$,
and these points satisfy the compatibility conditions of the
direct system). Hence to prove that $\phi_{\X}^{\rm int}$ induces
an isomorphism of the reduced schemes, it suffices to show that
$\underline{\phi_{\X}^{\rm int}}(R)$ is surjective for any reduced
$R$.

Note that $\phi_\X^\chi$ coincides with the restriction of
$\phi_V^\chi$ to $H_{\X, T}\subset H_{V, T}$  for any $\chi\in
\Sigma_V^{\rm int}\subset \Sigma_\X^{\rm int}$. By \cite[Theorem
5.6]{HS}, the map $\underline{\phi_{V}^{\rm int}}(R)$ is
surjective for any  reduced $R$, and it  follows that any element
$\{I_\chi\}_{\chi\in \Sigma^{\rm int}_\X}$ in  $\underline{
\X/\!\!_C} T(R)\subset \underline{ V/\!\!_C} T(R)$ gives an
element $\{I_\chi\}_{\chi\in \Sigma^{\rm int}_V}$ in
$\underline{H_{V, T}^{\rm int}}(R)$, i.e., $fI_{\chi_2}\subset
I_{\chi_1}$ for any $\chi_1, \chi_2 \in \Sigma_V^{\rm int}$ and
any $f\in S_{\chi_1-\chi_2}$. We have to prove that this condition
holds for any $\chi_1, \chi_2 \in \Sigma_\X^{\rm int}$. There
exists $c\in \N$ such that $c\chi_1, c\chi_2\in \Sigma^{\rm
int}_V$. For any $f'\in I_{\chi_2}$ we have $f^c(f')^c\in
I_{c\chi_1}$. By Lemma \ref{l1}, we see that the projection of
$\Sp ((R\otimes_k S^{(\chi_1)})/I^{(\chi_1)})$ to $\Sp R$ is a
locally trivial bundle with fiber $\A^1$. Consequently,
$(R\otimes_k S^{(\chi_1)})/I^{(\chi_1)}$ is reduced and $ff'\in
I^{(\chi_1)}$.
\end{proof}

\begin{remark} Note that restricting $\phi_{\X}$ to the main component $H_0$, we obtain a birational morphism of
toric $\T/T$-varieties from $H_0$ to $(\X/_{\!\! C}T)_0$.
\end{remark}

\begin{example} Let $V=\A^3$ where $\G^3$ and $T=\G$ act as in Example \ref{ex}, and let $\T=\G^2$ be embedded in $\G^3$ by $(t_1, t_2)\to (t_1, t_1, t_2)$. Consider the variety $\X=\overline{\T\cdot(1,1,1)}=\Sp S$,
where $S=k[x_1,x_2, x_3]/I_{\X}$ and $I_{\X}=(x_1-x_2)$.  So $H_{\X, T}$ is defined in $H_{\A^3, T}$ by the equation $z_1=z_2$.
We have the  homomorphisms  of  groups of characters
$$\Z^3=\Chi(\G^3)\stackrel{\pi'}{\to}\Z^2= \Chi(\T)\stackrel{\pi}{\to} \Z=\Chi(T)$$ and of monoids
$$\Omega_{\A^3}\to \Omega\to \Sigma,$$ where $\Omega_{\A^3}$ is the
monoid in $\Chi(\G^3)$ generated by characters with negative coordinates.

\begin{picture}(250, 70)
{\small \put(35, 50){$\Chi(\G^3)$}} \put(25, 20){\circle*{2}} \put(25, 20){\vector(-1, -1){15}} \put(0,
10){$\epsilon_1$} \put(25, 20){\vector(1, 0){20}} \put(40, 10){$\epsilon_2$} \put(10, 35){$\epsilon_3$} \put(25,
20){\vector(0, 1){20}} \put(70, 30){\vector(1, 0){15}} \put(73, 35){$\pi'$} \put(105, 20){\vector(0, 1){20}}
\put(105, 20){\circle*{2}} {\small \put(135, 50){$\Chi(\T)$} \put(110, 35){$e_2=\pi'(\epsilon_3)$} \put(120,
10){$e_1=\pi'(\epsilon_1)$} \put(135, 0){$=\pi'(\epsilon_2)$}} \put(235, 50){$\Chi(T)$} \put(225,
15){$\pi(e_1)$} \put(260, 15){$\pi(e_2)$} \put(105, 20){\vector(1, 0){20}} \put(180, 30){\vector(1, 0){15}}
\put(183, 35){$\pi$} \put(220, 25){\vector(1, 0){25}} \put(220, 25){\vector(1, 0){50}} \put(220,
25){\circle*{2}}

\end{picture}

Note that $\Sigma^{\rm int}_\X=\Sigma^{\rm int}_{\A^3}$ is the set of even numbers. The scheme $H^{\rm int}_{\A^3, T}$ is the closed subscheme in $\Pr^3$ defined by the equation $w_3^2=w_1w_2$, and $H^{\rm int}_{\X, T}$ is  defined by the equations $w_1=w_2=w_3$.
The isomorphism
$$\phi^{\rm int}_{\X}: H_{\X,T}^{\rm int}\to \X/_{\!\! C} T=\PP S$$ is the restriction of the isomorphism
$$\phi^{\rm int}_{\A^3}:
 H^{\rm int}_{\A^3, T}\to \A^3/_{\!\! C} T=\PP k[x_1,x_2,x_3],$$ where the inverse isomorphism is given by
$$(\phi^{\rm int}_{\A^3})^{-1}(x_1:x_2:x_3)=(x_1^2:x_2^2:x_1x_2:x_3).$$

Concerning the morphism $\phi_{\A^3}:
 H_{\A^3, T}\to \A^3/_{\!\! C}T$, note that $\phi_{\A^3}^{-1}(\X/_{\!\! C} T)$ is not contained in $H_{\X, T}$.
Indeed, consider the ideal $I=(x_1, x_2^2)\in \underline{H_{\A^3, T}}(k)$. We have
  $(I_{\X})_r\subset I_r$ for any even $r$, so $\phi_{\A^3}(I)\in \underline{\X/_{\!\! C}T}(k)$,
   but $I\notin \underline{H_{\X, T}}(k)$. \end{example}

\bigskip

\address{
Department of Higher Algebra \\
Faculty of Mechanics and Mathematics \\
Moscow State University \\
119992 Moscow\\
Russia
} {chuvasho@mccme.ru}

\end{document}